\documentclass [11pt]{amsart}

\usepackage {amsmath, amssymb, a4wide, epsfig, enumerate}
\usepackage[latin1]{inputenc}
\usepackage[frenchb]{babel}

\date{le 11 f{\'e}vrier 2005}
\thanks{{\it Codes Classif. AMS}: 53C12, 30C62}

\author{Romain Dujardin}
\title[Approximation sur les laminations]{Approximation des fonctions
  lisses sur certaines laminations}

\newcommand{\cc}{\mathbb{C}}

\newcommand{\dd}{\mathbb{D}}

\newcommand{\e}{\varepsilon}
\newcommand{\cv}{\rightarrow}

\newcommand{\fr}{\partial}
\newcommand{\om}{\Omega}
\newcommand{\set}[1]{\left\{#1\right\}}
\newcommand{\norm}[1]{\left\Vert#1\right\Vert}
\newcommand{\abs}[1]{\left\vert#1\right\vert}
\newcommand{\finc}{\subset \subset}
\newcommand{\cd}{\cc^2}
\newcommand{\pd}{{\mathbb{P}^2}}

\newcommand{\m}{{\bf M}}

\newcommand{\unsur}[1]{\frac{1}{#1}}
\newcommand{\el}{\mathcal{L}}

\newcommand{\rond}{\hspace{-.1em}\circ\hspace{-.1em}}

\newtheorem{theo}{Th{\'e}or{\`e}me} 
\newtheorem{prop} {Proposition} [section]
\newtheorem{thm}[prop] {Th{\'e}or{\`e}me} 
\newtheorem{cor}[prop]{Corollaire}
\newtheorem{lem}[prop] {Lemme}
\theoremstyle{definition}

\newtheorem{exam}[prop]{Exemple}

\newtheorem{rmk}[prop]{Remarque}

\begin{document}

\begin{abstract} Nous montrons que sur une lamination par surfaces de
  Riemann localement plong{\'e}e dans $\cd$, les fonctions $C^1$ sont localement
  limites uniformes de fonctions $C^1$ de l'espace ambiant, avec
  contr{\^o}le $L^p$ des d{\'e}riv{\'e}es le long des feuilles. On en d{\'e}duit qu'un
  courant domin{\'e} par un courant uniform{\'e}ment laminaire dans
  $\mathbb{C}^2$ est uniform{\'e}ment laminaire.
\end{abstract}

\maketitle

\section*{Introduction}
Les laminations par surfaces de Riemann sont des objets devenus usuels en
 dynamique holomorphe uni- et multi-dimensionnelle. En dimension complexe 2, toute
r{\'e}union de graphes disjoints dans un bidisque forme une telle lamination, souvent
appel{\'e}e  mouvement holomorphe. 

Il est bien connu que l'holonomie d'un mouvement holomorphe n'est en g{\'e}n{\'e}ral pas
lisse, ce qui est source de difficult{\'e}s. Le probl{\`e}me se pose en particulier pour la 
d{\'e}termination de la structure locale des ``courants structuraux'' (Sullivan \cite{sul}, 
voir en particulier la note p.236) ou ``courants dirig{\'e}s'' (Berndtsson-Sibony \cite{bs}
et Forn\ae ss-Sibony \cite{fs}) par une telle lamination.\\

Nous montrons ici le th{\'e}or{\`e}me de r{\'e}gularit{\'e} suivant pour une lamination par surfaces
de Riemann plong{\'e}e dans un ouvert $\om\subset\cd$ (Th{\'e}or{\`e}me
\ref{thm_approx}). On consultera  le paragraphe \ref{subs_C1} pour la notion de
fonction $C^1$ sur une lamination.  

\begin{theo}
Toute fonction de classe $C^1$ sur une telle lamination est limite
uniforme de restrictions de fonctions $C^1$  de l'espace ambiant, avec controle 
$L^p$ des d{\'e}riv{\'e}es le long des feuilles.
\end{theo}

Il est {\`a} noter que ce r{\'e}sultat ne r{\'e}sulte pas  directement
des th{\'e}or{\`e}mes usuels de r{\'e}gularit{\'e} des 
hom{\'e}omorphismes quasiconformes, qui ne traitent que des d{\'e}riv{\'e}es transverses aux 
feuilles (en particulier nous ne savons pas {\'e}tendre notre th{\'e}or{\`e}me au cas 
des laminations par hypersurfaces dans $\cc^n$). 

Le th{\'e}or{\`e}me a une  cons{\'e}quence  
sur la structure locale des courants laminaires (Th{\'e}\-or{\`e}me \ref{thm_dirig}) 
qui {\'e}tait {\`a} l'origine de ce travail. Nous obtenons en particulier le
r{\'e}sultat suivant (Corollaire \ref{cor_dirig}): 

\begin{theo} Si $T$ est un courant uniform{\'e}ment laminaire dans
  $\om\subset\cd$ et $S$
est un courant positif ferm{\'e} tel que $S\leq T$, alors $S$ est
uniform{\'e}ment laminaire.
\end{theo}

\medskip

Une motivation pour l'{\'e}tude de la r{\'e}gularit{\'e} des laminations 
par surfaces de Riemann en dimension  2, et des courants qu'elles supportent,
est la question suivante, qui demeure  toujours ouverte:

{\em Existe t'il une lamination par surfaces de Riemann dans 
$\pd$ munie d'une mesure transverse invariante?}

\medskip

La r{\'e}ponse (n{\'e}gative) {\`a} cette question a {\'e}t{\'e} donn{\'e}e 
 dans des cas particuliers :
\begin{itemize}
\item[-] par  Camacho-Lins Neto-Sad
\cite{cls} (voir {\'e}galement Deroin
\cite{de}) dans le cas o{\`u} la lamination est subordonn{\'e}e {\`a} un feuilletage holomorphe
singulier. Ces auteurs 
utilisent la construction de Bott d'une connexion plate dans le fibr{\'e}
 normal {\`a} la lamination,  qui doit pour ceci {\^e}tre de classe $C^1$.  
\item[-] Forn\ae ss-Sibony \cite{fs} d{\'e}montrent un r{\'e}sultat de non existence 
de telles laminations  en supposant soit une 
hypoth{\`e}se de r{\'e}gularit{\'e} transverse sur  la lamination (voir {\'e}galement
\cite{hm} dans ce cas), 
soit une hypoth{\`e}se 
de r{\'e}gularit{\'e} sur la mesure transverse. Ceci leur permet essentiellement 
de consid{\'e}rer le produit 
$T\wedge T$, o{\`u} $T$ est le cycle feuillet{\'e} associ{\'e}, pour parvenir 
{\`a} une contradiction. 
\end{itemize}

Dans les deux cas, 
c'est un probl{\`e}me de r{\'e}gularit{\'e} de l'holonomie qui emp{\^e}che
l'extension au cas g{\'e}n{\'e}ral. Les difficult{\'e}s soulev{\'e}es sont du m{\^e}me ordre que 
celles que nous rencontrerons dans cet article. 


\section{Mouvements holomorphes et courants}\label{sec_courants}

\subsection{Mouvements holomorphes} \label{subs_mouv}
Dans tout l'article, le cadre d'{\'e}tude sera le 
suivant : on se donne une famille $\el$ de graphes 
$L_\alpha =\set{w=\varphi_\alpha(z)}$ au dessus du 
disque unit{\'e} $\dd$ dans $\cd$, index{\'e}e par un param{\`e}tre transverse $\alpha\in\tau$, 
o{\`u} $\tau\subset\set{0}\times\dd$. On supposera de plus que pour tout 
$\alpha\in\tau$ et tout $z\in\dd$, $\abs{\varphi_{\alpha}(z)}\leq 1-\e<1$. 
En passant {\`a} la cl{\^o}ture de la famille de graphes, on pourra supposer que
$\tau$ est ferm{\'e}e. 

Dans un article classique, Ma{\~n}{\'e}, Sad et Sullivan \cite{mss}
ont montr{\'e} que l'application d'holonomie 
$$h_{z,z'}:(\set{z}\times\dd)\cap\el\cv(\set{z'}\times\dd)\cap\el$$ est continue et 
quasiconforme. Une telle application est en g{\'e}n{\'e}ral appel{\'e}e mouvement 
holomorphe. 
La famille de graphes $\el$ admet ainsi une structure de lamination
par surfaces de Riemann. Concluant une longue s{\'e}rie de travaux, 
S{\l}odkowski \cite{sl} a d{\'e}montr{\'e} qu'un mouvement holomorphe d{\'e}fini sur $\tau$ admet 
toujours une extension {\`a} un mouvement holomorphe de $\cc$ tout entier. La 
preuve {\`a} notre avis la plus accessible de ce th{\'e}or{\`e}me 
est celle de Chirka \cite{c}. 

Apr{\`e}s extension, l'holonomie est d{\'e}finie dans $\cc$ et quasiconforme, elle 
v{\'e}rifie donc une {\'e}quation de Beltrami
$$\frac{\fr h_{0,z}}{ \fr\bar w} = \mu^{z}(w) \frac{\fr h_{0,z}}{ \fr w},
\quad \abs{\mu^z(w)}<1.$$ Il est connu que $z\mapsto \mu^z(w)$ est 
holomorphe, et donc d'apr{\`e}s le lemme de Schwarz, 
$\abs{\mu^z}\leq \abs{z}$. En particulier le coefficient de Beltrami 
$\mu^z$ est arbitrairement petit si le temps de vie du mouvement holomorphe
est suffisamment court.

Une cons{\'e}quence est que l'application d'holonomie est H{\"o}ld{\'e}rienne. Il est possible de 
donner une preuve directe du caract{\`e}re H{\"o}lder de l'holonomie de la fa{\c c}on 
suivante : si $\alpha, \beta\in \tau$  sont distincts, la fonction 
$$k:s\mapsto -\log\unsur{2}\abs{\varphi_\alpha(s)-\varphi_\beta(s)}$$ est harmonique
et strictement positive. De l'in{\'e}galit{\'e} de Harnack
$$\frac{1-\abs{z}}{1+\abs{z}}\leq \frac{k(z)}{k(0)} \leq \frac{1+\abs{z}}{1-\abs{z}},$$ 
on d{\'e}duit que 
$$2^{\frac{-2\abs{z}}{1-\abs{z}}} \abs{\varphi_\alpha(0)-
\varphi_\beta(0)}^{\frac{1+\abs{z}}{1-\abs{z}}} \leq 
\abs{\varphi_\alpha(z)-\varphi_\beta(z)} \leq 
2^{\frac{2\abs{z}}{1+\abs{z}}} \abs{\varphi_\alpha(0)-
\varphi_\beta(0)}^{\frac{1-\abs{z}}{1+\abs{z}} } .$$ Un corollaire est
 que la constante de 
H{\"o}lder peut {\^e}tre choisie arbitrairement proche de 1, quitte {\`a} r{\'e}duire le disque 
de base. 

\subsection{Courants laminaires et dirig{\'e}s}
Si $\mu$ est une mesure positive finie sur $\tau$, la formule 
$$T=\int_\tau [L_\alpha]d\mu(\alpha)$$ d{\'e}finit un courant positif
ferm{\'e} de masse  localement finie 
dans $\dd^2$ qui est dit {\sl uniform{\'e}ment laminaire} et {\sl subordonn{\'e}} {\`a} la 
lamination $\el$. 

On peut d{\'e}finir une autre classe de courants associ{\'e}s {\`a} $\el$. Sullivan \cite{sul}
d{\'e}finit les 
{\sl courants structuraux} --que suivant Sibony, nous appellerons {\sl dirig{\'e}s}-- comme
le c{\^o}ne convexe ferm{\'e} engendr{\'e} par les courants de la forme
$i \delta_p (d\ell(p) \wedge d\overline\ell(p))$ o{\`u} $d\ell(p)$ 
est une forme $\cc$-lin{\'e}aire  dont le noyau est $T_pL$, o{\`u} $L$ est une feuille de $\el$ 
et $p\in L$. Dans la situation ci-dessus, on peut normaliser le choix de $d\ell$ en 
imposant que $d\ell(0,1)=1$; alors si $p\in L_\alpha$, 
$d\ell(p)=dw-\varphi'_\alpha(z)dz$. On v{\'e}rifie ais{\'e}ment que la
(1,0) forme $d\ell$ est continue, et un courant positif est dirig{\'e} par $\el$ ssi
$T\wedge d\ell=0$. Il est clair qu'un courant uniform{\'e}ment laminaire est dirig{\'e}. Plus
g{\'e}n{\'e}ralement, tout courant de la forme $fT$, o{\`u} $T$ est uniform{\'e}ment laminaire
et $f\in L^1_{loc}(\norm{T})$ ($\norm{T}$ d{\'e}signe la mesure trace de $T$)
est dirig{\'e}.

La question que nous nous posons est alors la suivante: est il vrai que tout
courant positif ferm{\'e} dirig{\'e} par un mouvement holomorphe $\el$ est 
uniform{\'e}ment laminaire? Si $\el$ est de classe $C^1$, le r{\'e}sultat se trouve
dans \cite{sul}. Dans le cas g{\'e}n{\'e}ral nous n'avons 
que le r{\'e}sultat partiel suivant.

\begin{thm}\label{thm_dirig}
Soit $T$ un courant uniform{\'e}ment laminaire subordonn{\'e} {\`a} une lamination 
par graphes au dessus de $\dd$. 
Alors tout courant positif ferm{\'e} de la forme $fT$, o{\`u} $f\in L^{1+\e}_{loc}(\norm{T})$,
$\e>0$,
est uniform{\'e}ment laminaire. 
\end{thm}

\begin{cor} \label{cor_dirig}
Si $T$ est uniform{\'e}ment laminaire et $S$ est un courant positif 
ferm{\'e} tel que $S\leq T$, alors $S$ est uniform{\'e}ment laminaire.
\end{cor}

\begin{proof} En effet d'apr{\`e}s le th{\'e}or{\`e}me de Radon-Nikodym 
il  existe une fonction mesurable $0\leq f\leq 1$ 
telle que $\norm{S}=f\norm{T}$. Par ailleurs le courant $T$ admet une d{\'e}composition
polaire 
$$T=\int \langle t(p), \cdot\rangle d\norm{T}(p), $$ o{\`u} $t(p)$ est un (1,1) vecteur
positif de norme 1, d{\'e}fini
$\norm{T}$ presque partout. $T$ {\'e}tant uniform{\'e}ment laminaire, le vecteur 
$t(p)$ est presque partout simple, i.e. de la forme $ie\wedge
\overline{e}$. On en d{\'e}duit 
ais{\'e}ment que $S= \int \langle t(p), \cdot\rangle fd\norm{T}(p)=fT$. 
\end{proof}

\begin{rmk} 
Le r{\'e}sultat {\'e}tant de nature locale, le th{\'e}or{\`e}me
 vaut {\'e}galement pour un courant subordonn{\'e} {\`a} 
une lamination plong{\'e}e dans une surface complexe. 
\end{rmk}

La preuve du th{\'e}or{\`e}me \ref{thm_dirig} que nous pr{\'e}sentons maintenant
est en quelque sorte une version rigoureuse du raisonnement
heuristique suivant
$$d(fT)=0 \Longrightarrow f \text{ constante le long des feuilles}.$$ Elle
repose de mani{\`e}re cruciale sur le th{\'e}or{\`e}me d'approximation
\ref{thm_approx}. 
\medskip

\noindent{\em Preuve du th{\'e}or{\`e}me.} On pose $S=fT$. Soit $\chi$ une
  fonction continue sur la transversale $\tau$. On peut prolonger
 $\chi$ en une 
fonction, toujours not{\'e}e $\chi$, constante le long des feuilles. 

\begin{lem}
Le courant $\chi S$ est  ferm{\'e}.
\end{lem}

\begin{proof} Soit $\phi$ une 1-forme test. On veut montrer que 
$\langle \chi S,d\phi \rangle=0$.  Supposons pour un instant que l'holonomie
soit lisse de classe $C^1$. Dans ce cas on peut {\'e}crire 
$$\langle \chi S,d\phi \rangle= -\langle S, d\chi\wedge \phi\rangle
=\int_\tau \left(\int_{L_\alpha} f d\chi\wedge \phi \right) d\mu(\alpha)=0$$ car
$\chi$ est $C^1$ et  constante le long des feuilles. 

Dans le cas g{\'e}n{\'e}ral  la fonction $\chi$ n'est plus $C^1$, et les calculs pr{\'e}c{\'e}dents
n'ont plus de sens. D'apr{\`e}s le th{\'e}or{\`e}me \ref{thm_approx}, pour tout $1\leq p <\infty$
$\chi$ est limite uniforme d'une suite de fonctions lisses $\chi_n$ de sorte que 
les d{\'e}riv{\'e}es de $\chi_n$ le long de $L_\alpha$
 tendent vers 0 dans $L^p$, uniform{\'e}ment en $\alpha$. On choisit $p>\frac{1+\e}{\e}$,
et  on {\'e}crit alors  
$$\langle \chi S,d\phi \rangle= \lim_{n\cv\infty} 
\langle \chi_n S,d\phi \rangle =
\int_\tau \left(\int_{L_\alpha} f d\chi_n\wedge \phi \right) d\mu(\alpha).$$ Ce dernier
terme tend vers 0 car $f\in L^{1+\e}$ et $p>\frac{1+\e}{\e}$. 
\end{proof}

Le th{\'e}or{\`e}me est alors une cons{\'e}quence directe de la caract{\'e}risation
suivante des courants uniform{\'e}ment laminaires.

\begin{lem}\label{lem_carac} 
Soit $S$ un courant positif ferm{\'e} port{\'e} par $\el$, tel que pour toute
fonction $\chi$ continue et constante le long des feuilles, $\chi S$
soit ferm{\'e}. Alors $S$ est uniform{\'e}ment laminaire.
\end{lem}

\begin{proof} On veut montrer que le courant
$S$ est dans le c{\^o}ne convexe ferm{\'e} engendr{\'e} par les $[L_\alpha]$. Supposons
le contraire. D'apr{\`e}s le th{\'e}or{\`e}me de Hahn-Banach, il existe une forme test 
$\phi$ telle que $\langle [L_\alpha], \phi \rangle<0$ pour tout $\alpha$ et 
$\langle S, \phi\rangle>0$. On consid{\`e}re alors
un recouvrement fini de $\tau$
par des ouverts $U_i$ de diam{\`e}tre inf{\'e}rieur {\`a} $\unsur{10}$, et une
partition 
continue de l'unit{\'e} 
$\theta_i$ associ{\'e}e. On prolonge les fonctions $\theta_i$ de 
mani{\`e}re  constante le long des feuilles. D'apr{\`e}s le lemme, chacun des courants
$\theta_iS$ est ferm{\'e}. De plus 
$$\langle S, \phi\rangle = \sum \langle\theta_i S, \phi\rangle > 0$$ donc 
il existe $i$ tel que $\langle\theta_i S, \phi\rangle > 0$. On pose $S_1=c_i\theta_iS$,
o{\`u} la constante $c_i>0$ est telle que $\m(S_1)=1$. Comme  $S_1=c_i\theta_i f T$, 
$S_1$ est de la forme $f_1 T$, o{\`u} $f_1\in L^{1+\e}_{loc}(\norm{T})$.

On refait la m{\^e}me construction pour $S_1$, avec des ouverts de diam{\`e}tre 
$<\unsur{100}$. En appliquant it{\'e}rativement le proc{\'e}d{\'e}, on obtient une suite de 
courants $S_k$ positifs ferm{\'e}s de masse 1, dont les supports convergent vers une 
certaine feuille $L$. Ainsi $S_k\cv c[L]$, avec $c>0$. Mais 
$\langle S_k,\phi\rangle>0$ alors que  $\langle [L], \phi \rangle<0$, ce qui est 
contradictoire.
\end{proof}

\begin{rmk}
Le lemme \ref{lem_carac} a pour cons{\'e}quence directe un
r{\'e}sultat d{\^u} {\`a} N. Sibony \cite{sibnonpub}: si $\el$ est une
lamination transversalement totalement discontinue, et $T$ est un
courant positif ferm{\'e} port{\'e} par $\el$, alors $T$ est uniform{\'e}ment
laminaire.
\end{rmk} 


\section{Le th{\'e}or{\`e}me d'approximation}\label{sec_approx}

\subsection{Structure $C^1$}\label{subs_C1}
Rappelons d'abord ce qu'on entend par fonction de classe $C^1$ sur une
lamination. Une fonction $f$ sur $\el$ sera dite
de classe $C^1(\el)$ si elle est
$C^1$ le long des feuilles et que les d{\'e}riv{\'e}es de $f$ le long des
feuilles sont transversalement continues. On norme l'espace $C^1(\el)$
de mani{\`e}re naturelle.

Il convient de remarquer que
m{\^e}me si $\el$ est un feuilletage de classe $C^1$, une fonction
$C^1(\el)$ n'est pas n{\'e}cessairement la restriction d'une fonction $C^1$ 
de l'espace ambiant. 

On a n{\'e}anmoins dans ce cas un r{\'e}sultat d'approximation 
\cite[Lemme 3.23]{de}, dont nous incluons une preuve succinte
pour la commodit{\'e} du lecteur. 

\begin{lem}\label{lem_C1}
Si $\el$ est de classe $C^1$, toute fonction $f\in C^1(\el)$ est
limite dans $C^1(\el)$ de restrictions de fonctions $C^1$ de l'espace
ambiant.
\end{lem}

\begin{proof} Apr{\`e}s redressement par un diff{\'e}omorphisme de classe
  $C^1$, on peut supposer que $\el$ est la lamination  naturelle  de 
$\dd\times\tau$. Si $g$ est une fonction continue sur $\tau$, on
dispose d'un prolongement continu $\mathcal{E}(g)$ de $g$ {\`a} $\cc$, o{\`u} 
$\mathcal{E}$ est l'op{\'e}rateur de prolongement de Whitney $C(\tau)\cv C(\cc)$
(voir \cite{st}). En particulier la fonction 
$$(z,\alpha) \mapsto
\mathcal{E}\left(f(z, \cdot)\right)(\alpha)$$ est d{\'e}finie dans
$\dd\times \cc$. Comme $\mathcal{E}$ est lin{\'e}aire et continu, le
prolongement est $C^1$ par rapport {\`a} la variable $z$. Il ne reste qu'{\`a}
convoler par un noyau r{\'e}gularisant de variable 
$\alpha$ pour obtenir l'approximation voulue.
\end{proof}

\begin{rmk} Par un argument analogue {\`a} la preuve du Th{\'e}or{\`e}me
  \ref{thm_dirig} on montre que si  $\el$  est une lamination 
par surfaces de Riemann plong{\'e}e satisfaisant la
  conclusion du lemme, tout courant positif ferm{\'e} dirig{\'e} par $\el$ est
  uniform{\'e}ment laminaire (voir {\'e}galement \cite[Proposition 3.22]{de}). 
\end{rmk}

\subsection{Approximation}
La notation $\el$ d{\'e}signe toujours une lamination par graphes au
dessus du disque unit{\'e}. On {\'e}crira $L\in\el$ pour dire que $L$ est une
feuille (plaque) de $\el$, et $\nabla_L$ d{\'e}signe le gradient le long
de $L$. 
Pour $p>1$ on d{\'e}signe par $W^{1,p}(\el)$ l'espace de
``Sobolev'' des fonctions continues sur $\el$ dont les d{\'e}riv{\'e}es dans
les feuilles sont dans $L^p$, et varient contin{\^u}ment au sens o{\`u}  
l'application $L \mapsto \nabla_L f$ est continue. On le munit de la
norme 
$$\norm{f}_{W^{1,p}} = \norm{f}_{L^\infty} + \sup_{L\in\el}
\norm{\nabla_Lf}_{L^p},$$ qui en fait un espace de Banach.

Le th{\'e}or{\`e}me  d'approximation est le suivant.

\begin{thm}\label{thm_approx}
Soit $1<p<\infty$ et $f\in C^1(\el)$. Il existe une suite  $(f_n)$
de fonctions lisses dans $\dd\times\cc$ telles que 
$$\norm{f-f_n}_{W^{1,p}}\underset{n\cv\infty}\longrightarrow 0$$ sur
tout compact.
\end{thm}

La preuve sera d{\'e}coup{\'e}e en {\'e}tapes et occupera le reste de la section.\\
  
{\em Localisation.} Remarquons d'abord qu'en recollant des approximations 
locales avec une partition de l'unit{\'e}, on peut se contenter de
r{\'e}soudre le probl{\`e}me dans un voisinage arbitraire $\set{0}\times D(0,r)$ de
$\set{0}\times\dd$. En particulier le coefficient de dilatation de
l'{\'e}quation de Beltrami (\ref{eq_beltrami}) ci-dessous
pourra {\^e}tre choisi arbitrairement petit. On
renormalise ensuite pour consid{\'e}rer {\`a} nouveau une lamination dans
$\dd\times\cc$. \\

{\em Extension.} 
Gr{\^a}ce au th{\'e}or{\`e}me de S{\l}odkowski, on {\'e}tend la
lamination, toujours not{\'e}e $\el$, 
{\`a} $\dd\times\cc$. On peut supposer que les feuilles sont
des droites horizontales hors du bidisque unit{\'e}. Comme au 
\S \ref{subs_mouv}, l'holonomie entre droites verticales 
v{\'e}rifie maintenant une {\'e}quation de Beltrami
\begin{equation}\label{eq_beltrami}
\frac{\fr h_{z,z'}}{ \fr\bar w} = \mu^{z,z'}(w) \frac{\fr h_{z,z'}}{ \fr w},
\quad \abs{\mu^{z,z'}(w)}\leq\kappa = \frac{2r}{1+r^2}<1,
\end{equation}
o{\`u} $r$ {\'e}t{\'e} d{\'e}fini ci-dessus. L'application d'holonomie 
$h^{z,z'}$ est  d{\'e}finie dans $\set{z}\times\cc$, et  par la suite
on {\'e}crira indistinctement  $h^{z,z'}(z,w)$ et $h^{z,z'}(w)$.
Remarquons que $\mu^{z,z'}$ est d{\'e}fini dans
$\cc$, et s'annule identiquement hors de $\dd$.\\ 

{\em Une premi{\`e}re approximation.} Sur $\el$ on dispose 
de deux syst{\`e}mes de
coordonn{\'e}es: les coordonn{\'e}es standard $(z,w)$ dans $\dd\times\cc$ et 
les coordonn{\'e}es $(z,\alpha)$, o{\`u} $\alpha\in \cc\simeq
\set{0}\times\cc$ indexe la feuille $L_\alpha$ {\`a} laquelle le point
courant $(z,w)$ appartient. 

L'application 
$$\Phi: (z,w)\longmapsto (z,\alpha)=(z,h^{z,0}(z,w))$$ 
est un hom{\'e}omorphisme qui redresse
$\el$, et qui est de classe $C^1(\el)$. La fonction $f\rond \Phi^{-1}$
est d{\'e}finie sur le ferm{\'e} $\dd\times\tau$ de la lamination redress{\'e}e. 
En utilisant le Lemme \ref{lem_C1} on approche $f\rond \Phi^{-1}$ par
des fonctions $g_n$, de classe 
$C^1$ (au sens usuel) dans $\dd\times \cc$ pour la topologie $C^1$ de la
lamination redress{\'e}e. On a maintenant $g_n\circ \Phi \cv f$ dans 
$W^{1,p}(\el)$, et les fonctions $g_n=f_n \rond \Phi^{-1}$ sont $C^1$ par
rapport {\`a} la 
variable $(z,\alpha)$ (mais pas par rapport {\`a} $(z,w)$). 
Quitte {\`a} remplacer $f$ par un $f_n$ assez
proche, il suffit donc de d{\'e}montrer le th{\'e}or{\`e}me quand $f\rond\Phi^{-1}
(z,\alpha)$ est   $C^1$ par rapport {\`a} $(z,\alpha)$. 
\footnote{Dans le cas,  
  suffisant pour le th{\'e}or{\`e}me \ref{thm_dirig}, 
o{\`u} $f$ est constante le long des feuilles, l'utilisation
  du Lemme \ref{lem_C1}  est superflue. Il suffit d'approcher uniform{\'e}ment
$f$ par des fonctions lisses dans
  la transversale $\set{0}\times \cc$ et d'{\'e}tendre les fonctions
  lisses en constantes le long des feuilles.}\\

{\em Convolution.} On convole le coefficient de Beltrami
$\mu^z(w)=\mu^{0,z}(w)$ par un noyau r{\'e}gulari\-sant $\theta_\e(w)$. Le
nouveau noyau $\mu_\e^z(w)=\mu^z(w)*\theta_\e(w)$ est toujours
holomorphe  en $z$ donc l'unique solution $h_\e^{0,z}(w)$ injective
et tangente {\`a}
l'identit{\'e} {\`a} l'infini de l'{\'e}quation de Beltrami correspondante est
l'holonomie d'une lamination par graphes
$\el_\e$ dans $\dd\times\cc$. D'apr{\`e}s un th{\'e}or{\`e}me classique (voir 
Ahlfors \cite{ah}),
$w \mapsto h_\e^{0,z}(w)$ est de classe $C^1$, et quand 
 $\e\cv 0$, $h_\e^{0,z}$ converge uniform{\'e}ment vers $h^{0,z}$,  
et les d{\'e}riv{\'e}es de $h_\e^{0,z}$ (par rapport {\`a} $w$) convergent en
norme $L^{p(\kappa)}$. L'exposant $p(\kappa)$ tend vers l'infini quand 
$\kappa \cv 0$. On choisit $\kappa$ de sorte que $p(\kappa)>p$, $p$
fix{\'e} par l'{\'e}nonc{\'e} du th{\'e}or{\`e}me. 

La famille de laminations $\el_\e$ tend vers $\el$ au sens o{\`u} la
feuille $L_\e(\alpha)$ issue de $(0, \alpha)$ converge vers
$L_\alpha$, uniform{\'e}ment relativement {\`a} $\alpha$. Soit $\Phi_\e$
l'hom{\'e}omorphisme de redressement de $\el_\e$ comme pr{\'e}c{\'e}demment.
On d{\'e}finit la fonction $f_\e$ par 
$$f_\e (z,w) = f\rond \Phi^{-1}\rond \Phi_\e (z,w) \Leftrightarrow
f_\e \rond \Phi_\e^{-1}(z,\alpha) =  f\rond \Phi^{-1} (z,\alpha),$$
autrement dit apr{\`e}s redressement, ``$f_\e=f$''. Il est clair que 
$f_\e$ tend uniform{\'e}ment vers $f$ quand $\e\cv 0$. Nous allons montrer que 
cette approximation r{\'e}sout notre probl{\`e}me.\\

{\em $f_\e$ est $C^1$.} Il faut prendre garde au fait
qu'une lamination dont l'holonomie est $C^1$ n'est pas un 
feuilletage $C^1$ en g{\'e}n{\'e}ral (voir le tr{\`e}s {\'e}clairant 
\S 6 de \cite{psw}). Ici, nous allons
v{\'e}rifier que l'application de carte
$$\Phi_\e^{-1}:(z,\alpha)\mapsto (z,w)=(z,( h_\e^{0,z})(\alpha))
$$ est de classe $C^1$ par rapport {\`a} la variable $(z,\alpha)$. Puisque 
$h_\e^{0,z}(\alpha)=\varphi^\e_\alpha(z)$ est holomorphe en $z$, il suffit
 de v{\'e}rifier que les d{\'e}riv{\'e}es dans la direction verticale  de 
$h_\e^{0,z}(\alpha)$ sont continues par rapport {\`a} $(z,\alpha)$. 

L'argument est le suivant (voir \cite[pp. 22-24]{cg} ou \cite{ah}): 
soit $\mu^z_\e(\alpha)$ le 
coefficient de Beltrami associ{\'e} {\`a} $h_\e^{0,z}(\alpha)$; $\mu^z_\e(\alpha)$ est 
$C^1$ en $\alpha$ et holomorphe en $z$. En particulier 
  $\frac{\fr}{\fr \alpha} \mu^z_\e(\alpha)$ est continue en $(z,\alpha)$. 
  On peut {\'e}crire  
$\frac{\fr}{\fr \alpha} h_\e^{0,z}(\alpha) = e^\varphi$, o{\`u} $\varphi$
v{\'e}rifie une {\'e}quation de type Beltrami de variable $z$, dont le second membre est 
$\frac{\fr}{\fr \alpha} \mu^z_\e(\alpha)$. Ceci implique que les d{\'e}riv{\'e}es 
de $h_\e^{0,z}$ dans la direction verticale 
sont uniform{\'e}ment born{\'e}es et H{\"o}lderiennes d'exposant $1-\frac{2}{p}$ en 
$\alpha$. Par ailleurs $\frac{\fr}{\fr \alpha} h_\e^{0,z}$, vue comme
fonction de $\alpha$, varie contin{\^u}ment dans 
$L^p$ relativement {\`a} la variable $z$, et donc contin{\^u}ment pour la topologie 
de la convergence uniforme, par {\'e}quicontinuit{\'e}. 
On en d{\'e}duit la continuit{\'e} de $\frac{\fr}{\fr \alpha} h_\e^{0,z}$ en la
variable $(z,\alpha)$.\\
 
{\em Convergence en norme $W^{1,p}$.} On fixe maintenant une feuille 
$L$ de $\el$, et on va montrer que $\norm{\nabla_L f_\e-\nabla_L f}_{W^{1,p}}
\cv 0$. Le lecteur v{\'e}rifiera ais{\'e}ment que l'uniformit{\'e} en $L$ est
cons{\'e}quence de la preuve. 

Quitte {\`a} redresser par une application holomorphe, on peut toujours supposer
que $L=L_0=(w=0)$. On pose $\widetilde{f_\e}= f_\e\circ \Phi_\e^{-1}$ (et de m{\^e}me pour $f$)
de sorte que $\widetilde{f_\e}(z,\alpha)= f_\e (z,w)$. Par d{\'e}finition de 
$f_\e$, on a $$f_\e(z, 0) = \widetilde{f_\e} (z, h^{z,0}_\e(z,0)) =
\widetilde{f}(z, h^{z,0}_\e(z,0)),$$ car $h^{z,0}_\e(z,0)$ est le projet{\'e} du point $(z,0)$
sur $\set{0}\times \cc$ le long de $\el_\e$. Cette projection
n'est certainement pas injective, puisque $\el_\e$ est essentiellement parall{\`e}le
{\`a} $\el$.  

Deux cas se pr{\'e}sentent: soit $L_0$ est {\'e}galement une feuille de $\el_\e$, auquel 
cas $f_\e(z,0)= f(z,0)$ et il n'y a rien {\`a} d{\'e}montrer; dans le cas contraire, les 
intersections de $L$ et $\el_\e$ sont transverses sauf en au plus un nombre localement 
fini de points \cite[Lemme 6.4]{bls}, donc en un point g{\'e}n{\'e}rique de $L$, 
$L$ est une transversale holomorphe locale {\`a} $\el_\e$.

Le point crucial de la preuve est le lemme suivant. Il est peut {\^e}tre d{\'e}j{\`a} connu 
des sp{\'e}cialistes, mais nous ne l'avons pas trouv{\'e} dans la litt{\'e}rature. Rappelons que 
pour toute paire de 
 droites verticales $\set{z}\times\cc$ et $\set{z'}\times \cc$, l'holonomie
$h^{z,z'}_\e$ est quasiconforme et de dilatation major{\'e}e par $\kappa$.

\begin{lem}\label{lem_dilat}
Soient $D_1$ et $D_2$ deux transversales holomorphes locales {\`a} une 
feuille $L_\e^0$ de $\el_\e$. Alors l'holonomie
$$h^{1,2}: D_1\subset U_1\longrightarrow U_2\subset D_2$$ est quasiconforme
de dilatation major{\'e}e par $\kappa$.
\end{lem}

En d'autres termes, la dilatation ne d{\'e}pend pas des transversales, 
pourvu qu'elles soient holomorphes.  

\begin{proof}[Preuve du lemme]  
La lamination $\el_\e$ est en fait un feuilletage de 
classe $C^1$, donc $h^{1,2}$ est diff{\'e}rentiable.  Pour simplifier les notations, 
dans la suite de la preuve nous  omettrons les indices $\e$.
Si $L$ est une feuille
proche de $L^0$ on pose  pour 
$i=1,2$, $p_i=(z_i,w_i)=D_i\cap L$, de sorte que $h^{1,2}(p_1)=p_2$. Les droites
$\set{z_1}\times\cc$ et $\set{z_2}\times \cc$ sont des transversales respectives
en $p_1$ et $p_2$, donc au voisinage de $p_1$, on peut d{\'e}composer 
$h^{1,2}$ en $h^{z_2, 2}\rond h^{z_1, z_2} \rond h^{1,z_1}$ o{\`u} $h^{1,z_1}$ 
(resp. $h^{z_2, 2}$) envoie
localement $D_1$ sur $\set{z_1}\times \cc$ (resp.  $\set{z_2}\times \cc$ sur $D_2$).

L'observation\footnote{Ceci
reste vrai m{\^e}me si $\el$ n'est pas lisse.} 
 est que $h^{1,z_1}$ est diff{\'e}rentiable au point 
$p_1$, et sa diff{\'e}rentielle est la projection de $T_{p_1}D_1$ sur 
$\set{z_1}\times\cc$ parall{\`e}lement {\`a} $T_{p_1}L$. C'est
 une application $\cc$-lin{\'e}aire, et donc son coefficient de Beltrami au point 
 $p_1$ est nul. Il en va de m{\^e}me pour $h^{z_2, 2}$ en $p_2$.
 Ainsi le coefficient de dilatation de $h^{1,2}$ en $p_1$ est major{\'e} par $\kappa$.
 \end{proof}

Nous pouvons maintenant terminer la preuve du th{\'e}or{\`e}me. On veut majorer
$\widetilde{f}(z,0)-\widetilde{f}(z, h^{z,0}_\e(z,0))$ en norme $W^{1,p}$. 
D'apr{\`e}s la premi{\`e}re approximation, $\widetilde{f}$ est $C^1$, donc il suffit de montrer
que les d{\'e}riv{\'e}es de $\pi_\e:z\mapsto h^{z,0}_\e(z,0)$ sont petites en norme $L^p$. 

D'apr{\`e}s le lemme, au voisinage de tout 
point de transversalit{\'e} de $L=\set{w=0}$ et $\el_\e$, la projection 
$\pi_\e$ est un hom{\'e}omorphisme quasiconforme de coefficient de
dilatation $\leq \kappa$. Donc hors d'un ensemble discret de points de
$L$, $\pi_\e$ v{\'e}rifie une {\'e}quation de 
Beltrami 
\begin{equation}\label{eq_beltrami2}
\frac{\fr \pi_\e}{\fr \bar z} = \nu \frac{\fr \pi_\e}{\fr  z}, 
\text{ avec } \abs{\nu}\leq \kappa .
\end{equation}
On identifie $L$ et $\dd$ et on {\'e}tend $\nu$ {\`a} $\cc$ par z{\'e}ro. 
Toutes les solutions de cette {\'e}quation sur $\dd$
sont de la forme $\pi_\e= u_\e\rond k_\e$ o{\`u} $k_\e$ est l'unique 
hom{\'e}omorphisme tangent {\`a} l'infini solution  de (\ref{eq_beltrami2}) et 
$u_\e$ est holomorphe sur $k_\e(\dd)$ \cite{do}.  En particulier 
$\norm{\nabla k_\e}_{L^p}\leq C$ o{\`u} $C$ ne d{\'e}pend que de $\kappa$, et 
les modules de continuit{\'e} (H{\"o}lder) $k_\e$ et $k_\e^{-1}$ ne d{\'e}pendent
{\'e}galement que de $\kappa$. 

Par ailleurs, $\el_\e$ est proche de $\el$, donc les feuilles de $\el_\e$ issues
de $L$ coupent $\set{0}\times \cc$ pr{\`e}s de 0. En particulier l'image de $\pi_\e$, 
c'est {\`a} dire celle  de $u_\e$ est de petit diam{\`e}tre. On fixe un disque
l{\'e}g{\`e}rement plus petit
$D(0,1-\delta)\finc \dd$. On a alors l'estimation
 $$\mathrm{dist}(\fr (k_\e(D(0,1-\delta)), \fr (k_\e(\dd)))\geq
 c(\kappa) \delta^{\theta(\kappa)},$$
 et donc d'apr{\`e}s la formule de Cauchy, les d{\'e}riv{\'e}es de $u_\e$ sont uniform{\'e}ment
 petites sur $k_\e(D(0,1-\delta))$. On en conclut que les d{\'e}riv{\'e}es de 
 $k_\e=  u_\e\rond k_\e$ sont petites en norme $L^p$ sur tout
 compact. \hfill $\square$\\


\bigskip

\noindent{\sc \small UFR de math{\'e}matiques,
Universit{\'e} Paris 7,
Case 7012,
2 place Jussieu,
75251 Paris cedex 05, France.}\\
{\tt\small dujardin@math.jussieu.fr}

\end{document}